\documentclass[11pt,a4paper]{article}
\usepackage[utf8]{inputenc}
\usepackage{graphicx,verbatim,array,multicol,courier}
\usepackage{amsmath,amssymb,amsthm,mathrsfs, mathtools}
\usepackage{color, graphics, graphicx, siunitx}
\usepackage{amsfonts, dsfont, bm}
\usepackage{natbib} 
\usepackage{a rydshln}
\usepackage{geometry}
\usepackage{enumitem,float}
\usepackage{lscape}
\usepackage[pdftex,bookmarks,colorlinks]{hyperref}
\usepackage[dvipsnames]{xcolor}
\usepackage{tabularx}
\usepackage{siunitx}
\usepackage[labelformat=simple]{subfig}
\usepackage[linesnumbered,ruled,vlined]{algorithm2e}
\usepackage{booktabs}

\newlist{inparaenum}{enumerate}{2}
\setlist[inparaenum,1]{label=(\alph*)}
\setlist[inparaenum,2]{label=(\roman{inparaenumi}\emph{\alph*})}

\makeatletter
\def\adl@drawiv#1#2#3{%
        \hskip.5\tabcolsep
        \xleaders#3{#2.5\@tempdimb #1{1}#2.5\@tempdimb}%
                #2\z@ plus1fil minus1fil\relax
        \hskip.5\tabcolsep}
\newcommand{\cdashlinelr}[1]{%
  \noalign{\vskip\aboverulesep
           \global\let\@dashdrawstore\adl@draw
           \global\let\adl@draw\adl@drawiv}
  \cdashline{#1}
  \noalign{\global\let\adl@draw\@dashdrawstore
           \vskip\belowrulesep}}
\makeatother

\newcommand{\blind}{0}

\numberwithin{equation}{section}
\theoremstyle{definition}
\newtheorem{defi}{Definition}[section]
\newtheorem{cond}[defi]{Condition}
\theoremstyle{plain}
\newtheorem{theo}[defi]{Theorem}
\newtheorem{prop}[defi]{Proposition}
\newtheorem{lem}[defi]{Lemma}
\newtheorem{cor}[defi]{Corollary}

\theoremstyle{remark}

\theoremstyle{example}
\newtheorem{ex}[defi]{Example}

\newcommand{\diff}{\mathrm{d}}

\definecolor{navy}{rgb}{0,0,0.502}
\definecolor{brown}{rgb}{0.59, 0.29, 0.0}
\def\indic{\mathds{1}}

\hypersetup{colorlinks,%
	citecolor=blue,%
	filecolor=green,%
	linkcolor=red,%
	urlcolor=violet,%
}

\newcommand{\hell}{{\mathscr{H}}}
\newcommand{\kulb}{{\mathscr{K}}}
\newcommand{\variational}{{\mathscr{V}}}

\newcommand{\pdiv}{{\mathscr{D}}}
\newcommand{\gdist}{{\mathscr{D}}}
\newcommand{\TEmeasure}{{\mathcal{P}}}
\newcommand{\GPmeasure}{{\mathscr{P}}}

\newcommand{\Real}{\mathbb{R}}
\newcommand{\DoA}{\mathcal{D}}
\newcommand{\Prob}{\mathbb{P}}
\newcommand{\Borel}{\mathbb{B}}

\title{Strong Convergence of Peaks Over a Threshold}
\author{S. A. Padoan\\
	Department of Decision Sciences, Bocconi University, Italy\\
	and\\ 
	S. Rizzelli\\
	Department of Statistical Sciences, Catholic University, Italy
}

\begin{document}

\maketitle
\begin{abstract}
Extreme Value Theory plays an important role to provide approximation results for the extremes of 
a sequence of independent 
random  variables when their distribution is unknown. An important one is given by the {generalised Pareto distribution} $H_\gamma(x)$ as an approximation of the distribution $F_t(s(t)x)$ of the excesses over a threshold $t$, where $s(t)$ is a suitable norming function. In this paper we study the rate of convergence of $F_t(s(t)\cdot)$ to $H_\gamma$ in variational and Hellinger distances and translate it into that regarding the Kullback-Leibler divergence between the respective densities. 
%
\end{abstract}

\\


%
\section{Introduction}
%
	
Extreme value theory (EVT)  develops probabilistic models and methods for describing the random behaviour of extreme observations that rarely occur. These theoretical foundations are very important for studying practical problems in environmental, climate, insurance and financial fields (e.g., \citealp{embrechts2013modelling}; \citealp{dey2016extreme}; \citealp{kulik2020}), to name a few. 

In the univariate setting, the most popular approaches for statistical analysis are the so-called {block maxima} (BM) and {peaks over threshold} (POT) \citep[see e.g.][for a review]{axel21}. Let $X_1,\ldots,X_n$ be independent and identically distributed (i.i.d.) random variables according to a common distribution $F$. The first approach concerns the modelling of $k$ sample maxima derived over  blocks of a certain size $m$, i.e. $M_{m,i}=\max(X_{(i-1)m+1},\ldots,X_{im}), i\in\{1,\ldots,k\}$. In this case, under some regularity conditions \citep[e.g.][Ch. 1]{dehaan+f06}, the weak limit theory establishes that $F^m(a_mx+b_m)$ converges pointwise to $G_\gamma(x)$ as $m\to\infty$, for every continuity point $x$ of $G_\gamma$, where $ G_\gamma$ is the {generalised extreme value} (GEV) distribution, $a_m>0$ and $b_m$ are suitable norming constants for each $m=1,2,\ldots$ and $\gamma\in\Real$ is the so-called tail index, which describes the tail heaviness of $F$ \citep[e.g.][Ch. 1]{dehaan+f06}. The second method concerns the modelling of $k$ random variables out of the $n$ available that exceed a high threshold $t$, or, equivalently, of $k$ threshold excesses $Y_j$, $j=1,\ldots,k$, which are i.i.d. copies of $Y=X-t|X>t$. In this context, the {generalised Pareto} (GP) distribution, say $H_\gamma$, appears as weak limit law of appropriately normalised high threshold exceedances, i.e. for all $x>0$, $F_t(s(t)x)$ converges pointwise to $H_{\gamma}(x)$ as $t\to x^*$, for all the continuity points $x$ of $H_\gamma(x)$, where $F_t(x)=\Prob(Y\leq x)$ and $s(t)>0$ is a suitable scaling function for any $t\leq x^*$, with $x^*={\sup}(x: F(x)<\infty)$. This result motivates the POT approach, which was introduced decades ago by the seminal paper  \citet{balkema1974residual}. Since then, few other convergence results emerged. For instance, the uniform convergence of $F_t(s(t)\, \cdot \,)$ to $H_\gamma$ and the coresponding convergence rate have been derived  by \citet{pickands1975statistical} and \citet{raoult03}, respectively. Similar results but in  Wasserstein distance have been recently established  by \citet{bobbia2021coupling}. As for the GEV distribution, more results are available. In particular, there are sufficient conditions to ensure, in addition to weak convergence, that  $F^m(a_m\, \cdot \, +b_m)$ converges to $G_\gamma$ for example uniformly and in variational distance and the density of $F^m(a_m\, \cdot \, +b_m)$ converges pointwise, locally uniformly and uniformly to that of $G_\gamma$ (e.g. \citealp{falk2010}, Ch. 2; \citealp{resnick2007}, Ch. 2). 

The main contribution of this article is to provide new convergence results that can be useful in practical problems for the POT approach.
Motivated by the utility in the statistical field to asses the asymptotic accuracy of estimation procedures, we study stronger forms of convergence than the pointwise one, as $\lim_{t\to x^*}\gdist(F_t(s(t)\, \cdot \,),  H_{\gamma})=0$,  where $\gdist(\,\cdot \,;\, \cdot \,)$ is either the variational distance, the Hellinger distance or the Kullback-Leibler divergence.
%
%
In particular, we provide upper bounds for the rate of convergence to zero of $\gdist(F_t(s(t)\,\cdot \,);  H_{\gamma})$ in the case that $\gdist(\,\cdot \,;\, \cdot \,)$ is the variational and Hellinger distance, and
further translate them into bounds  on Kullback-Leibler divergence between the densities of $F_t(s(t)\cdot)$ and $H_{\gamma}$, respectively.
We also pinpoint cases where reentering of exceedances is necessary to reach the optimum rate, namely where $F_t(s(t)\cdot+c(t))$ has to be considered in place of $F_t(s(t)\cdot)$, for a suitable real valued function $c(t)$.

\if1\blind
{
Estimators of the tail index $\gamma$ (and other related quantities) are typically defined as functionals of the random variables $(Y_1,\ldots,Y_k)$, as for instance the popular Hill  \citep{hill1975simple}, Moment \citep{dekkers1989moment}, Pickands \citep{pickands1975statistical}, Maximum Likelihood (ML, \citealp{jenkinson1969statistics}), Generalised Probability Weighted Moment (GPWM, \citealp{hosking1985estimation}) estimators, to name a few. 
In real applications, the distribution $F$ is typically unknown and so is $F(s(t)\, \cdot \,)$.
Although, for large $t$, $H_\gamma$ provides a model approximation for $F_t(s(t)\,\cdot \,)$, when one wants to derive asymptotic properties as the consistency and especially the rate of convergence of the tail index estimators (or other related quantities), still the fact that (after rescaling) the random variables $(Y_1,\ldots,Y_k)$ are actually distributed according to $F_t(s(t)\,\cdot \,)$ needs to be taken into account, which makes asymptotic derivations quite burdensome.
These are even more complicated  if $t$ is determined on the basis of the $(k+1)$-th largest order statistic of the original sample $X_1,\ldots, X_n$, which is the most common situation in practical applications. In this case, the threshold is in fact random and, up to rescaling, {\color{magenta}}$F_t(s(t)\, \cdot \,)$ only gives a conditional model for the variables $Y_j$ given a fixed value $t$ of the chosen statistic.
Asymptotic properties for POT methods have been studied in the last fifty years, see for example \citet{hall1984best},  \citet{drees1998optimal}, \citet{dekkers1993optimal} and the reference{} therein.  

Leveraging on our strong convergence results we can show that, for random sequences (such as sequences of estimators) convergence results in probability that hold under the limit model $H_\gamma$, are also valid for a rescaled sample of excesses  over a large order statistic. Precisely, we show that the distribution of the latter, up to rescaling and reordering, is {\it contiguous} to that of an ordered i.i.d. sample from $H_\gamma$ \citep[e.g.,][Ch. 6.2]{vdv2000}. 
%
As a by product of this result, one can derive the consistency and rate of convergence  of a tail index estimator (or an estimator of a related quantity) by defining it as a functional of the random sequence $(Z_1,\dots,Z_k)$ which is distributed according the limit model $H_\gamma$, and, if density of $F_t(s(t)\cdot)$ satisfies some regularity conditions, then the same asymptotic results hold even when such estimator is defined through the sequence of excesses.
This approach simplifies a lot the computations as asymptotic properties are easily derivable under the limit model. 
}\fi

The article is organised as follows, Section \ref{sec:background} of the paper provides a brief summary of the probabilistic context on which our results are based. Section \ref{sec:strong_results} provides our new results on strong convergence to a Pareto model. 
Section  \ref{sec:proofs} provides the proofs of the main results.

%
\section{Background}\label{sec:background}
%
	
Let $X$ be a random variable with a distribution function $F$ that is in the domain of attraction of the GEV distribution $G_\gamma$, shortly denoted as $F\in\DoA(G_\gamma)$. This means that there are norming constants $a_m>0$ and $b_m\in\Real$ for $m=1,2,\ldots$ such that 
\begin{equation}\label{eq:GEV}
\lim_{m\to\infty}F^m(a_m x + b_m)=\exp\left(-\left(1+\gamma x\right)^{-1/\gamma}\right)=:G_\gamma(x),
\end{equation}
for all $x\in\Real$ such that $1+\gamma x>0$, where $\gamma\in\Real$, and this is true if only if there is a scaling function $s(t)>0$ with $t< x^*$ such that
\begin{equation}\label{eq:GP}
\lim_{t\to x^{*}}F_t(s(t)x)=1-\left(1+\gamma x\right)^{-1/\gamma}=:H_\gamma(x),
\end{equation}
e.g., \citet[][Theorem 1.1.6]{dehaan+f06}. The densities of $H_\gamma$ and $G_\gamma$ are 
	$$
	h_\gamma(x)=\left(1+\gamma x \right)^{-(1/\gamma+1)}
	$$
	and
	$$
	g_{\gamma}(x)=G_\gamma(x)h_\gamma(x),
	$$
respectively. Let $U(v):=F^{\leftarrow}(1-1/v)$, for $v\geq 1$, where $F^{\leftarrow}$ is the left-continuous inverse function of $F$ and $G^{\leftarrow}(\exp(-1/x))=(x^\gamma-1)/\gamma$. Then, we recall that the first-order condition in formula \eqref{eq:GEV} is equivalent to the limit result
	\begin{equation}\label{eq:quant}
			\lim_{v\to\infty}\frac{U(vx)-U(v)}{a(v)}=\frac{x^\gamma-1}{\gamma},
	\end{equation}
	for all $x>0$, where $a(v)>0$ is a suitable scaling function. In particular, we have that $s(t)= a(1/(1-F(t)))$, see \citet[][Ch. 1]{dehaan+f06} 
 	for possible selections of the function $a$.
	
	A stronger convergence form than that in formula \eqref{eq:GP} is the uniform one, i.e.
	$$
	\sup_{x\in[0,\frac{x^*-t}{s(t)})}\left|F_t(s(t)x)-H_\gamma(x)\right|\to 0, \quad t\to x^*.
	$$
	In case of distributions $F$ with finite end-point $x^*$, the following slightly more general form of convergence is also of interest
	$$
	\sup_{x\in[0,\frac{x^*-t-c(t)}{s(t)})}\left|F_t(s(t)x+c(t))-H_\gamma(x)\right|\to 0, \quad t\to x^*,
	$$
	for a centering function $c(t)$ satisfying $c(t)/s(t)\to0$ as $t\to x^*$.
	To establish the speed at which $F_t(s(t) x)$ or $F_t(s(t)x+c(t))$ converges uniformly to $H_\gamma(x)$, \citet{raoult03} 
	relied on a specific formulation of the well-known  second-order condition. In its general form, the second order condition requires the existence of a positive function $a$ and a positive or negative function $A$, named rate function, such that $\lim_{v\to \infty}|A(v)|=0$ and
	$$
	\lim_{v\to\infty}\frac{\frac{U(vx)-U(v)}{a(v)}-\frac{x^\gamma-1}{\gamma}}{A(v)}=D(x), \quad x>0,
	$$
	where $D$ is a non-null function  which is not a multiple of $(x^\gamma-1)/\gamma$,  see \citet[][Definition 2.3.1]{dehaan+f06}. The rate function $A$ is necessarily regularly varying at infinity with index $\rho\leq0$, named second-order parameter \citep[][Theorem 2.3.3]{dehaan+f06}.
	In the sequel, we use the same specific form of second order condition of \citet{raoult03} to obtain decay rates for stronger metrics than uniform distance between distribution functions.

	%
	\section{Strong results for POT}\label{sec:strong_results}
	%
	%
	In this section, we discuss strong forms of convergence for the distribution of 
	renormalised
	exceedances over a threshold.
	First, in Section \ref{sec:total_variation}, we discuss convergence to a GP distribution in variational and Hellinger distance, drawing a connection with known results for density convergence of normalized maxima.
	In Section \ref{sec:hellinger_distance} we quantify the speed of convergence in variational and Hellinger distance.
	%
	Moreover, we show how these can be used to also bound Kullback-Leibler divergences.
	Throughout, for a twice differentiable function $W(x)$ on $\mathbb{R}$, we denote with $W'(x)=(\partial/\partial x)W(x)$ and $W''(x)=(\partial^2/\partial x^2)W(x)$ the first and second order derivatives, respectively.
\subsection{Strong convergence under classical assumptions}\label{sec:total_variation}
Let the distribution function $F$ be twice differentiable. In the sequel, we denote  $f=F'$, $g_m=(F^m(a_m\,\cdot \,+b_m))'$ and $f_t=F_t'$.
Under the following classical von Mises-type conditions
	\begin{eqnarray}\label{eq:vm}
		\nonumber
		\lim_{x\to\infty}\frac{xf(x)}{1-F(x)}=\frac{1}{\gamma}, \quad \gamma>0,\\
		\lim_{x\to x^*}\frac{(x^*-x)f(x)}{1-F(x)}=-\frac{1}{\gamma}, \quad \gamma<0,\\
		\nonumber
		\lim_{x\to x^*}\frac{f(x)\int_x^{x^*}(1-F(v)dv)}{(1-F(x))^2}=0,\quad \gamma=0,
	\end{eqnarray}
	we know that the first-order condition in formula \eqref{eq:quant} is satisfied and it holds that  
	\begin{equation}\label{eq:conv}
		\lim_{v\to\infty}	v a(v)f(a(v)x+U(v))=(1+\gamma x)^{-1/\gamma-1}
	\end{equation}
 	locally uniformly for  $(1+\gamma x) >0$. Since the equality $
 	g_m(x)=F^{m-1}(a_mx+b_m) h_m(x)$ holds true, with $b_m=U(m)$, $a_m=a(m)$ and $h_m(x)=ma_mf(a_m x +b_m),
 	$ 
 	and since $F^{m-1}(a_mx+b_m)$ converges to $G_\gamma(x)$ locally uniformly as $m\to\infty$,
 	the convergence result in formula \eqref{eq:conv} thus implies that $g_m(x)$ converges to $g_\gamma(x)$ locally uniformly
 	\citep[][Ch. 2.2]{resnick2007}.
		
On the other hand, the density pertaining to $F_t(s(t)x)$ is
$$
l_t(x):=f_t(s(t)x)s(t)=\frac{s(t)f(s(t)x+t)}{1-F(t)}
$$
and, setting  $v=1/(1-F(t))$, we have $a(v)=s(t)$ and $v\to \infty$ as $t\to x^*$.	
Therefore, a further implication of the convergence result in formula \eqref{eq:conv} is that $l_t(x)$ converges to $h_\gamma(x)$ locally uniformly for $x>0$, if $\gamma \geq 0$, or $x\in (0,-1/\gamma)$, if $\gamma <0$.
	%
In turn,  by Scheffe's lemma we have
$$
\lim_{t\to x^*}\variational(\TEmeasure_t, \GPmeasure)=0,
$$
where 
$$
\variational(\TEmeasure_t, \GPmeasure)=\sup_{B\in\Borel}\left|\TEmeasure_t(B) - \GPmeasure(B)\right|
$$
is the total variation distance between the probability measures 
$$
\TEmeasure_t(B):=\Prob\left(\frac{X-t}{s(t)}\in B \bigg{|} X>t\right)\; \text{ and } \;\GPmeasure(B):=\Prob(Z\in B),
$$
and where $Z$ is a random variable with distribution $H_\gamma$ and $B$ is a set in the Borel $\sigma$-field of $\Real$, denoted by $\Borel$. 
Let
$$
\hell^2(l_t; h_\gamma):=\int \left[
\sqrt{l_t(x)}-\sqrt{h_\gamma(x)}
\right]^2 \diff x
$$
be the square of the Hellinger distance. It is well {known} that the Hellinger and total variation distances are related as
\begin{equation}\label{eq:hell_var_relation}
\hell^2(l_t; h_\gamma)\leq 2\variational(\TEmeasure_t, \GPmeasure)\leq 2\hell(l_t; h_\gamma) ,
\end{equation}
see e.g. \citet[][Appendix B]{ghosal2017}.
Therefore, the  conditions in formula \eqref{eq:vm} ultimately entail that also the Hellinger distance between the density of rescaled peaks over a threshold $l_t$ and the GP density $h_\gamma$ converges to zero as $t \to x^*$. In the next subsection we introduce a stronger assumption, allowing us to also quantify the speed of such convergence.
	%
	
%
%
%
\subsection{Convergence rates}\label{sec:hellinger_distance}
As in \cite{raoult03} we rely on the following assumption, in order to derive the convergence rate for the variational and Hellinger distance.	
%
\begin{cond}\label{cond:SO}
Assume that $F$ is twice differentiable. Moreover, assume that there exists $\rho\leq 0$ such that
$$
A(v):= \frac{vU''(v)}{U'(v)}+1-\gamma
$$
defines a function of constant sign near infinity, whose absolute value $|A(v)|$ is regularly varying as  $v\to\infty$ with index of variation  $\rho$.
\end{cond}
When Condition \ref{cond:SO} holds then the classical von-Mises conditions in formula \eqref{eq:vm} are also satisfied for the cases where $\gamma$ is positive, negative or equal to zero, respectively. Furthermore,  Condition \ref{cond:SO} implies that an appropriate scaling function for the exceedances of a high threshold $t<x^*$, which complies with the equivalent first-order condition \eqref{eq:GP}, is defined as 
$$
s(t)=(1-F(t))/f(t).
$$
With such a choice of the scaling function $s$, we establish the following results.
\begin{theo}\label{theo:hellrate}
Assume Condition \ref{cond:SO} is satisfied.
Then, there exist constants 
$c>0$, $\alpha_j>0$ with $j=1,2$, $K>0$ and $t_0 < x^*$, depending on $\gamma$, such that 
%
\begin{equation}\label{eq:hellbound}
\frac{\hell^2(
	\tilde{l}_t
	; h_\gamma)}{K |A(v)|^2} \leq S(v)
\end{equation}
for all $t\geq t_0$, where $v=1/(1-F(t))$,  
$$
\tilde{l}_t =\begin{cases}
	l_t, \hspace{7.9em} \gamma\geq 0\\
	l_t\left(\, \cdot \, +
	\frac{x^*-t}{s(t)}+\frac{1}{\gamma}
	\right), \quad \gamma <0
\end{cases}
$$
 and
$$
S(v)=
1-|A(v)|^{\alpha_1 }+4 \exp\left(
c |A(v)|^{\alpha_2}
\right).
$$

\end{theo}
Note that $\tilde{l}_t$ is the density of $F_t(s(t) \cdot +c(t))$, with centering function
\begin{equation}\label{eq:center}
c(t):=\begin{cases}
	0, \hspace{7.4em} \gamma \geq 0\\
	x^*-t +\gamma^{-1}s(t), \quad \gamma <0
\end{cases},
\end{equation}
for $t<x^*$.
Given the relationship between the total variation and Hellinger distances in \eqref{eq:hell_var_relation}, with obvious adaptations when a nonnul recentering is considered, the following result is a direct consequence of Theorem \ref{theo:hellrate}.
\begin{cor}\label{cor:totalvar_rate}
Under the assumptions of Theorem \ref{theo:hellrate},
for all $t\geq t_0$
$$
\variational(
\widetilde{\TEmeasure}_t
, \GPmeasure)\leq  |A(v)| \sqrt{K S(v)}{,}
$$
with $\widetilde{\mathcal{P}}_t$ the probability measure pertaining to $\tilde{l}_t$.
\end{cor}
Theorem \ref{theo:hellrate}  implies that {when $\gamma \geq 0$} the Hellinger and variational distances of the probability density and measure of rescaled exceedances from their GP distribution counterparts are bounded from above by $C  |A(v)|$, for a positive constant $C$, as the threshold $t$ approaches the end-point $x^*$.
Since for a fixed $x \in \cap_{t\geq t_0}(0,\frac{x^*-t}{s(t)})$ it holds that
$$
\left|F_t(s(t)x)-H_\gamma(x)\right| \leq \variational(\TEmeasure_t, \GPmeasure)
$$
and since \citet[][Theorem 2(i)]{raoult03} implies that  $|F_t(s(t)x)-H_\gamma(x)|/|A(v)|$ converges to a positive constant, there also exists $c>0$ such that, for all large $t$, $c|A(v)|$ is a lower bound for variational and Hellinger distances. Therefore, since 
$$
c|A(v)|\leq \variational(\TEmeasure_t, \GPmeasure)\leq  \hell(l_t; h_\gamma) \leq C|A(v)|,
$$
the decay rate of variational and Hellinger distances is precisely $|A(v)|$ as $t\to x^*$.
When $\gamma <0$, analogous considerations apply to $\tilde{l}_t$ and $\widetilde{\TEmeasure}_t$.  With the following results, we give precise indications of when a recentered version of $l_t$  is necessary to achieve the optimal rate.
\begin{prop}\label{prop:boundGPD}
Under the assumptions of Theorem \ref{theo:hellrate}, when $\gamma <0$ there are constants $c_j$, $j=3,4$, and $t_1 < x^*$, depending on $\gamma$, such that	for all $t >t_1$
$$
c_3 |A(v)|^{-1/2\gamma} < \hell
\left(
h_\gamma, h_\gamma\left(
\, \cdot \, -\mu_t 
\right)
\right)
<  c_4 |A(v)|^{\min(1, -1/2\gamma)},
$$
where $\mu_t:=c(t)/s(t)$ and $c(t)$ is as in the second line of \eqref{eq:center}.
\end{prop}
\begin{cor}\label{cor:rate}
Under the assumptions of Theorem \ref{theo:hellrate}: 
\begin{inparaenum}
	\item when $-1/2 \leq \gamma <0$, there are constants $c_5>0$ and $t_2 < x^*$, depending on $\gamma$, such that for all $t>t_2$
	$$
	\hell\left(
	l_t, h_\gamma
	\right) \leq c_5 |A(v)|;
	$$
	\item \label{cor:rate2}
	when $ \gamma <-1/2$, there are constants $c_6>0$ and $t_3 < x^*$, depending on $\gamma$, such that for all $t>t_3$
	$$
	\hell\left(
	l_t, h_\gamma
	\right) \geq c_5 |A(v)|^{-1/2\gamma}.
	$$
\end{inparaenum}
\end{cor}
According to Corollary \ref{cor:rate}\ref{cor:rate2}, the density $l_t$  of rescaled exceedances $Y/s(t)$ does not achieve the optimal convergence rate $|A(V)|$ whenever $\gamma_0 <-1/2$, in which case the rate is only of order $|A(v)|^{-1/2\gamma}$. In simple terms, this is due to the fact that, 
when $\gamma $ is negative, the supports of $l_t$ and $h_\gamma$ can be different and the approximation error is affected by the amount of probability mass in the unshared region of points. We recall indeed that
the end-point $(x^*-t)/{s}(t)$ of $l_t$ converges to $-1/{\gamma}$ as $t$ approaches $x^*$ at rate $A(v)$ \citep[e.g.,][Lemma 4.5.4]{dehaan+f06}. Nevertheless, for $t<x^*$ it can be that $(x^*-t)/{s}(t)>-1/\gamma$ or $(x^*-t)/{s}(t)<-1/\gamma $. 
In turn, when $\gamma$ is smaller than -1/2, the approximation error due to support mismatch has a dominant effect.
However, if scaled exceedances are shifted by subtracting the quantity  $\mu_t$, in this case the upper end-point of the density $\tilde{l}_t$ is the same of that of $h_\gamma$, hence no support mismatch occurs and the optimal convergence rate is achieved also in the case where $\gamma <-1/2$.

		%
		

A further implication of Theorem \ref{theo:hellrate} concerns the speed of convergence to zero of the Kullback-Leibler divergence 
	$$
	\kulb(\tilde{l}_t; h_\gamma):= \int \ln \left\lbrace\tilde{l}_t(x)/h_\gamma(x)\right\rbrace\tilde{l}_t(x) \diff x,
	$$
	and the divergences of higher order $p\geq 2$
	$$
	\pdiv_p(\tilde{l}_t; h_\gamma):= \int \left|\ln \left\lbrace\tilde{l}_t(x)/h_\gamma(x)\right\rbrace \right|^p\tilde{l}_t(x) \diff x{.}
	$$
%
	%
	%
%
{U}sing the uniform bound on density ratio provided in Lemma \ref{lem:boundens}
we are able to translate the upper bounds on the squared Hellinger distance $\hell^2({\tilde{l}_t},h_\gamma)$ into upper bounds on the Kullback-Leibler divergence $\kulb(\tilde{l}_t; h_\gamma)$ and higher order divergences $\pdiv_p(\tilde{l}_t; h_\gamma)$.
\begin{cor}\label{cor:kulbrate}
Under the assumptions of Theorem \ref{theo:hellrate} with in particular $\rho <0$ and $\gamma \neq 0$, 
there exist  constants {$M>0$ and  $t_4<x^*$, depending on $\gamma$,  such that for all $t\geq t_4$}	
	\begin{inparaenum}
		\item 
		$
		\kulb(\tilde{l}_t; h_\gamma)\leq 2M {K S(v) |A(v)|^2}
		$
		\item 
		$
		\pdiv_p(\tilde{l}_t; h_\gamma)\leq 2 p! {K S(v) |A(v)|^2}
		$
		with 	$p\geq 2$.
	\end{inparaenum}

	\end{cor}
To extend the general results in Lemma \ref{lem:boundens} and Corollary \ref{cor:kulbrate} to the case of $\gamma=0$ seems to be technically over complicated. 
Nevertheless, there are specific examples where the properties listed in such lemmas are satisfied, such as the following one.
	\begin{ex}
		Let $F(x)=\exp(-\exp(-x))$, $x\in \mathbb{R}$, be the Gumbel distribution function. In this case, Condition \ref{cond:SO} is satisfied with $\gamma=0$ and $\rho=-1$, so that Theorem \ref{theo:hellrate} applies to this example, and for an arbitrarily small $\epsilon>0$ we have
		$$
		l_t(x)/h_0(x)\leq \exp(\exp(-t))<1+\epsilon
		$$
		for all $x>0$ and suitably large $t$. Hence, the bounded density ratio property is satisfied and it is still possible to conclude that $\pdiv_p({l}_t; h_0)/|A(v)|^2$ and $\kulb({l}_t; h_0)/|A(v)|^2$ can be bounded from above as in Corollary \ref{cor:kulbrate}.
	\end{ex}
\section{Proofs}\label{sec:proofs}
\subsection{Additional notation}
For $y >0$, we denote $T(y)=U(e^y)$ and, for $t < x^*$, we define the functions  
$$
p_t(y)=
\begin{cases}
\frac{T(y+T^{-1}(t))-t}{s(t)}-\frac{e^{\gamma y}-1}{\gamma}, \hspace{4.5em} {\gamma > 0}\\
\frac{T(y+T^{-1}(t))-t}{s(t)}-y,
\hspace{6.3em}  \gamma=0\\
{\frac{T(y+T^{-1}(t))-x^* -\gamma^{-1}s(t)}{s(t)}-\frac{e^{\gamma y}-1}{\gamma},} \quad {\gamma < 0}
\end{cases},
$$
with $s(t)=(1-F(t))/f(t)$,
and
$$
q_t(y)=\begin{cases}
\frac{1}{\gamma}\ln \left[
1+\gamma e^{-\gamma y} p_t(y)
\right],\quad \gamma \neq 0\\
p_t(y),\hspace{7.9em} \gamma=0
\end{cases}.
$$
{Moreover, for $t<x^*$ we set
$$
\tilde{x}_t^*=\begin{cases}
	\frac{x^*-t}{\tilde{s}(t)}, \quad \gamma \geq 0\\
	-\frac{1}{\gamma}, \hspace{1.6em} \gamma <0
\end{cases}.
$$
Fuerthermore,}
for $x \in (0, x^*-t)$, we let
$
\phi_t(x)=T^{-1}(x+t)-T^{-1}(t).
$
Finally, for $x \in \mathbb{R}$, $\gamma \in \mathbb{R}$, $\rho \leq 0${, we set
$$
I_{\gamma,\rho}(x)=
\begin{cases}
	\int_0^x e^{\gamma s}\int_0^s e^{\rho z} \diff z \diff s,  \hspace{2.2em} \gamma \geq 0\\
	-\int_x^\infty e^{\gamma s}  \int_0^s e^{\rho z} \diff z \diff s, \quad \gamma <0 .
\end{cases}
$$
}
%
%
\subsection{Auxiliary results}
In this section we provide some results which are auxiliary to the proofs of the main ones, presented in Section \ref{sec:strong_results}. Throughout, for Lemmas \ref{lem:aux1}--\ref{lem:densexpand2}, Condition \ref{cond:SO} is implicitly assumed to hold true. The proofs are provided in the supplementary material document.
{ In particular, Lemmas \ref{lem:aux1}--\ref{lem:aux2}, provided next, are directly used in the proof of our main result, Theorem \ref{theo:hellrate}. 
}
\begin{lem}\label{lem:aux1}
	{For every $\varepsilon>0$ and every $\alpha>0
	$,
	there exist $x_1 < x^*$ and $\kappa_1>0$ (depending on $\gamma$) such that} for all $t \geq x_1$ and $y \in (0,-\alpha \ln 	\vert
	A(e^{T^{-1}(t)})
	\vert )$
		$$
			e^{ -\kappa_1
			\vert
			A(e^{T^{-1}(t)})
			\vert 
			e^{2\varepsilon y}
		   }
		<
		e^{{-}q_t(y)} <
		e^{ \kappa_1
			\vert
			A(e^{T^{-1}(t)})
			\vert 
			e^{2\varepsilon y}}.
		$$
\end{lem}
\if1\blind
{	
\begin{proof}
	By Lemma 5 in \cite{raoult03}, for all $\varepsilon>0$ there exists $x_0$ such that for all $t\in(x_0,x^*)$ and $y>0$,
\begin{equation*}
	\begin{split}
		e^{-\gamma x}\vert p_t(y)\vert\leq (1+\varepsilon)\vert
		A(e^{T^{-1}(t)})
		\vert I_{\gamma,\rho}(y)e^{(\gamma -\varepsilon)y}.
	\end{split}
\end{equation*}
Moreover, for a positive constant $\vartheta_1$
\begin{equation*}
		I_{\gamma,\rho}(y) e^{(\gamma -\varepsilon)y}\leq 
			\vartheta_1 e^{2\varepsilon y}.
\end{equation*}
Combining these two inequalities, we deduce that
\begin{equation}\label{eq:main1}
		e^{-\gamma y}\vert p_t(y)\vert\leq 	
		(1+\varepsilon)\vert
		A(e^{T^{-1}(t)})
		\vert \vartheta_1
		e^{2\varepsilon y}.
\end{equation}
As a consequence,
for any $\alpha>0$ there exists a constant $\vartheta_2$ such that
\begin{equation}\label{eq:main2}
	\sup_{y \in (0,-\alpha \ln 	\vert
		A(e^{T^{-1}(t)})
		\vert )}
	e^{-\gamma y}\vert p_t(y)\vert\leq \vartheta_2
	\vert
	A(e^{T^{-1}(t)})
	\vert^{1-2\varepsilon\alpha}.
\end{equation}
Therefore, choosing $\varepsilon$  sufficiently small, $e^{-\gamma y}\vert p_t(y)\vert$ converges to zero uniformly over the interval $(0,-\alpha \ln 	\vert
A(e^{T^{-1}(t)})
\vert )$ as $t\to x^*$.

It now follows that, if $y \in (0,-\alpha \ln 	\vert
A(e^{T^{-1}(t)})
\vert )$ and  $t>x_1$ for a sufficiently large value $x_1<x^*$, when $\gamma\neq 0$ a first-order Taylor expansion of the logarithm at $1$ yields
\begin{equation*}
	\begin{split}
		\vert
		q_t(y)
		\vert&=
		\left\vert
		\frac{1}{\gamma}
		\frac{\gamma e^{-\gamma y}p_t(y)}{1+\vartheta(t,y)\gamma e^{-\gamma y}p_t(y)}
		\right\vert\\
		&\leq 
		\vartheta_4 \vert
		A(e^{T^{-1}(t)})
		\vert e^{2\varepsilon y},
\end{split}
\end{equation*}
where $\vartheta(t,y)\in(0,1)$ and $\vartheta_4$ is a positive constant, while when $\gamma=0$ it holds that
\begin{equation*}
\begin{split}
	\vert
	q_t(y)
	\vert&=
	e^{\gamma y}	e^{-\gamma y}\vert
	p_t(y)
	\vert\\
	&\leq \vartheta_5 \vert
	A(e^{T^{-1}(t)})
	\vert e^{2\varepsilon y},
\end{split}
\end{equation*}
where  $\vartheta_5$ is a positive constant. The  result in the statement is a direct consequence of the last two inequalities.
%
\end{proof}
}\fi
\begin{lem}\label{lem:aux2}
	{For every $\varepsilon>0$ and every $\alpha>0$, 
	there exist $x_2 < x^*$ and $\kappa_2>0$ (depending on $\gamma$) such that} for all $t \geq x_2$ and $y \in (0,-\alpha \ln 	\vert
	A(e^{T^{-1}(t)})
	\vert )$
		$$
		e^{ -\kappa_2
			\vert
			A(e^{T^{-1}(t)})
			\vert 
			e^{2\varepsilon y}
		}
		<
		1+q_t'(y) <
		e^{ \kappa_2
			\vert
			A(e^{T^{-1}(t)})
			\vert 
			e^{2\varepsilon y}
		}.
		$$
\end{lem}
\if1\blind
{
\begin{proof}
	If $\gamma\neq 0$
	$$
	1+q_t'(y)=\frac{\exp\left\lbrace
		\int_{e^{T^{-1}(t)}}^{e^{y+T^{-1}(t)}}
		\frac{A(u)}{u}\diff u
		\right\rbrace
	}{1+\gamma e^{-\gamma y}p_t(y)},
	$$
	while if $\gamma=0$
	$$
	1+q_t'(y)={\exp\left\lbrace
		\int_{e^{T^{-1}(t)}}^{e^{y+T^{-1}(t)}}
		\frac{A(u)}{u}\diff u
		\right\rbrace
	}.
	$$
	Therefore, if $y \in (0,-\alpha \ln 	\vert
	A(e^{T^{-1}(t)})
	\vert )$ and  $t>x_2$ for a sufficiently large value $x_2<x^*$,  using the bounds in formulas
	\eqref{eq:main1}--\eqref{eq:main3} and choosing a suitably small $\varepsilon$ we deduce
	\begin{equation*}
	\begin{split}
	1+q_t'(y)&\leq 
	\frac{\exp\left\lbrace
		\int_{e^{T^{-1}(t)}}^{e^{y+T^{-1}(t)}}
		\frac{A(u)}{u}\diff u
		\right\rbrace
	}{1-\indic(\gamma\neq 0) \vert \gamma \vert e^{-\gamma y}\vert p_t(y)\vert}\\
	&\leq \exp\left\lbrace
	y\vert
	A(e^{T^{-1}(t)})
	\vert
	\right\rbrace \times
	\begin{cases}
	\frac{1}{1-\omega_1 \vert
		A(e^{T^{-1}(t)})
		\vert  e^{2\varepsilon y}}, \hspace{2em} \gamma \geq 0\\
	\frac{1}{1-\omega_2 \vert
		A(e^{T^{-1}(t)})
		\vert  e^{(\gamma -\varepsilon)y}}, \quad \gamma <0
	\end{cases}\\
	&\leq 
	\begin{cases}
	\exp\left\lbrace
	\omega_3 \vert
	A(e^{T^{-1}(t)})
	\vert  e^{2\varepsilon y}
	\right\rbrace, \hspace{2.2em} \gamma \geq 0\\
	\exp\left\lbrace
	\omega_4 \vert
	A(e^{T^{-1}(t)})
	\vert  e^{(\gamma-\varepsilon)y}
	\right\rbrace, \quad \gamma <0
	\end{cases}
	\end{split}
	\end{equation*}
	for positive constants $\omega_i$, $i=1,\ldots,4$. Similarly,
	\begin{equation*}
	\begin{split}
	1+q_t'(y)&\geq 
	\frac{\exp\left\lbrace
		\int_{e^{T^{-1}(t)}}^{e^{y+T^{-1}(t)}}
		\frac{A(u)}{u}\diff u
		\right\rbrace
	}{1+\indic(\gamma\neq 0) \vert \gamma \vert e^{-\gamma y}\vert p_t(y)\vert}\\
	&\geq \exp\left\lbrace
	-y\vert
	A(e^{T^{-1}(t)})
	\vert
	\right\rbrace \times
	\begin{cases}
	\frac{1}{1+\omega_5 \vert
		A(e^{T^{-1}(t)})
		\vert  e^{2\varepsilon y}}, \hspace{2em} \gamma \geq 0\\
	\frac{1}{1+\omega_6 \vert
		A(e^{T^{-1}(t)})
		\vert  e^{(\gamma -\varepsilon)y}}, \quad \gamma <0
	\end{cases}\\
	&\geq 
	\begin{cases}
	\exp\left\lbrace-
	\omega_7 \vert
	A(e^{T^{-1}(t)})
	\vert  e^{2\varepsilon y}
	\right\rbrace, \hspace{2.2em} \gamma \geq 0\\
	\exp\left\lbrace
	-\omega_8 \vert
	A(e^{T^{-1}(t)})
	\vert  e^{(\gamma-\varepsilon)y}
	\right\rbrace, \quad \gamma <0
	\end{cases}
	\end{split}
	\end{equation*}
	for positive constants $\omega_i$, $i=5,\ldots,8$. The result now follows.
\end{proof}
}\fi
{Lemmas \ref{lem:bigO}--\ref{lem:densexpand2} are auxiliary to the proof of Lemma \ref{lem:boundens}, which in turn is the cornerstone of the proof of Corollary \ref{cor:kulbrate}. Their statemens are however reported here below, as they can be of interest on their own.}
\begin{lem}\label{lem:bigO}
	If $\gamma>0$ and $\rho<0$, there exists a regularly varying function $\mathcal{R}$ with negative index $\varrho$ such that, defining the function 
	$$
	\eta(t):=\frac{(1+\gamma t)f(t)}{1-F(t)} -1, 
	$$
	as $v\to\infty$, 
	$
	\eta(U(v))=O(\mathcal{R}(v)).
	$
\end{lem}
\if1\blind
{
\begin{proof}
	Let $v_0>0$ satisfy $U(v_0)\neq 0$ and $U'(v_0)\neq 0$. Then, for $v >v_0$ it holds that \begin{equation*}
	\begin{split}
	\eta(U(v))&=\frac{1+\gamma U(v)}{vU'(v)}-1\\
	&=\frac{1+\gamma U(v_0)}{vU'(v)}
	+\gamma \int_{v_0}^v\frac{U'(r)}{vU'(v)}\diff r-1.
	\end{split}
	\end{equation*}
	Moreover, by definition of $A$, we have the identity
	\begin{equation*}
	\begin{split}
	\gamma \int_{v_0}^v\frac{U'(r)}{vU'(v)}\diff r-1&=\gamma\int_{v_0/v}^1\frac{U'(zv)}{U'(v)}\diff z-1\\
	&=\int_{v_0/v}^{1} \gamma z^{\gamma-1}\left[
	\exp\left\lbrace
	-\int_z^1 \frac{A(vu)}{u}\diff u
	\right\rbrace
	-1
	\right]\diff z-\left(\frac{v_0}{v}\right)^\gamma.
	\end{split}
	\end{equation*}
	Therefore, denoting by $\mathcal{R}_2(v)$ the first term on the right-hand side and setting
	$$
	\mathcal{R}_1(v)=\frac{1+\gamma U(v_0)}{vU'(v)}-\left(\frac{v_0}{v}\right)^\gamma,
	$$
	we have $\eta(U(v))=\mathcal{R}_1(v)+\mathcal{R}_2(v)$. On one hand, the function $\mathcal{R}_1(v)$ is regularly varying of order $-\gamma$. On the other hand, for any $\beta \in(0,1)$, the function $\mathcal{R}_2(v)$ can be decomposed as follows
	\begin{equation*}
	\begin{split}
	\mathcal{R}_2(v)&=\int_{v_0/v}^{v^{-(1-\beta)}}+\int_{v^{-(1-\beta)}}^1
	\gamma z^{\gamma-1}\left[
	\exp\left\lbrace
	-\int_z^1 \frac{A(vu)}{u}\diff u
	\right\rbrace
	-1
	\right]\diff z\\
	&=:\mathcal{R}_{2,1}(v)+\mathcal{R}_{2,2}(v).
	\end{split}
	\end{equation*}
	Assuming that $A$ is ultimately positive and selecting $v_0$ suitably large, we have
	\begin{equation*}
	\begin{split}
	\vert 
	\mathcal{R}_{2,1}(v)
	\vert&\leq \int_{v_0/v}^{v^{-(1-\beta)}}
	\gamma z^{\gamma-1}\left[1-
	\exp\left\lbrace
	-\frac{A(vz)}{z}
	\right\rbrace
	\right]\diff z\\
	&=O(v^{-\gamma(1-\beta)})
	\end{split}
	\end{equation*}
	and 
	\begin{equation*}
	\begin{split}
	\vert 
	\mathcal{R}_{2,2}(v)
	\vert&\leq \int_{v^{-(1-\beta)}}^1
	\gamma z^{\gamma-1}\left[1
	-z^{A(v^\beta)}
	\right]\diff z\\
	&=O(v^{-\gamma(1-\beta)}\vee A(v^\beta)).
	\end{split}
	\end{equation*}
	Consequently, there exists a regularly varying function $\mathcal{R}$ of index $\varrho= \gamma(\beta-1) \vee \rho\beta$ complying with the property in the statement as $v\to\infty$.

	Similarly, if $A$ is ultimately negative, choosing $\beta$ such that $\beta <2\gamma$ and $v_0$ suitably large, we have
	\begin{equation*}
	\begin{split}
	\vert 
	\mathcal{R}_{2,1}(v)
	\vert&\leq \int_{v_0/v}^{v^{-(1-\beta)}}
	\gamma z^{\gamma-1}\left[
	u^{A(v_0)}
	-1
	\right]\diff z \\
	&=O(v^{-(\gamma -\beta/2)(1-\beta)})
	\end{split}
	\end{equation*}
	and
	\begin{equation*}
	\begin{split}
	\vert 
	\mathcal{R}_{2,2}(v)
	\vert&\leq \int_{v^{-(1-\beta)}}^1
	\gamma z^{\gamma-1}\left[
	z^{A(v^\beta)}-1
	\right]\diff z\\
	&=O(v^{-(\gamma -\beta/2)(1-\beta)} \vee \vert A(v^\beta)\vert )
	\end{split}
	\end{equation*}
	as $v\to\infty$. Hence, there exists a regularly varying function $\mathcal{R}$ of index $\varrho= (\beta-1)(\gamma-\beta/2) \vee \rho\beta$ complying with the property in the statement. The proof is now complete.
\end{proof}
}\fi
\begin{lem}\label{lem:densexpand1}
	If $\gamma>0$ and $\rho<0$, there exists $x_3\in(0,\infty)$ and $\delta >0$ such that, for all $x\geq x_3$,
	$$
	f(x)=   h_\gamma(x)\left[
	1+O(\{1-H_\gamma(x)\}^\delta)
	\right].
	$$
\end{lem}
\if1\blind
{
\begin{proof}
	Let $\mathcal{R}^*(t):=\mathcal{R}(1/(1-F(t)))$,  where $\mathcal{R}$ is as in Lemma \ref{lem:bigO}. Then $\mathcal{R}^*(t)$ is regularly varying of index $\varrho/\gamma$ \citep[][Proposition 0.8(iv)]{resnick2007}. In turn, by Karamata's theorem \citep[e.g,][Proposition 0.6(a)]{resnick2007} we have that for a large $t^*$
	$$
	\int_{ t^* }^\infty \frac{\vert \eta(t)\vert}{1+\gamma t} \diff t<\infty
	$$
	and thus, by Proposition 2.1.4 in \cite{falk2010}, we conclude that
	\begin{equation}\label{eq:tailequiv}
	\tau:=\lim_{t\to\infty}\frac{1-F(t)}{1-H_\gamma(t)} \in (0,\infty).
	\end{equation}
	As a consequence, for any $\delta \in (0,-\varrho)$, as $t\to \infty$
	\begin{equation*}
	\begin{split}
	\mathcal{R}^*(t)&\sim \mathcal{R}\left(
	\frac{1}{\tau(1-H_\gamma(t))}
	\right)\\
	&=O(\{1-H_\gamma(t)\}^{\delta}).
	\end{split}
	\end{equation*}
	The conclusion now follows by Proposition 2.1.5 in \cite{falk2010}.
\end{proof}
}\fi
\begin{lem}\label{lem:bigO2}
	If $\gamma<0$ and $\rho<0$, there exists a a regularly varying function $\tilde{\mathcal{R}}$ with negative index $\tilde{\varrho}=(-1)\vee(-\rho/\gamma) $ such that, defining the function 
	$$
	\tilde{\eta}(y):=\frac{(1-\gamma y)f(x^*-1/y)}{[1-F(x^*-1/y)]y^2} -1, 
	$$
	as $y\to\infty$, 
	$
	\tilde{\eta}(y)=O(\tilde{\mathcal{R}}(y)).
	$
\end{lem}
\if1\blind
{
\begin{proof}
	By definition,
	\begin{equation*}
	\begin{split}
	\tilde{\eta}\left(y\right)&= \frac{f(x^*-1/y)}{[1-F(x^*-1/y)]y^2}-\gamma\left[
	\frac{f(x^*-1/y)}{y(1-F(x^*-1/y))}+\frac{1}{\gamma}
	\right]\\
	&=:\tilde{\eta}_1\left(y\right)
	+\tilde{\eta}_2\left(y\right).
	\end{split}
	\end{equation*}
	On one hand, we have that, as $y\to\infty$ 
	$$
	\tilde{\eta}_1\left(y\right)=O(1/y).
	$$
	On the other hand, for $v>1$ we have the identity
	$$
	\tilde{\eta}_2\left(\frac{1}{x^*-U(v)}\right)=\int_1^\infty \gamma z^{\gamma-1}
	\left[
	1-\exp\left\lbrace
	\int_1^z \frac{A(uv)}{u}\diff u
	\right\rbrace
	\right]\diff z.
	$$
	Hence, if $A$ is ultimately positive,
	\begin{equation*}
	\begin{split}
	\tilde{\eta}_2\left(\frac{1}{x^*-U(v)}\right) &\leq -\gamma \int_1^\infty z^{\gamma -1}(z^{A(v)}-1)\diff z\\
	&=O(A(v))
	\end{split}
	\end{equation*}
	while, if $A$ is ultimately negative,
	\begin{equation*}
	\begin{split}
	\left\vert \tilde{\eta}_2\left(\frac{1}{x^*-U(v)}\right) \right\vert &\leq \gamma A(v) \int_1^\infty z^{\gamma -1}\ln z\diff z\\
	&=O(\vert A(v) \vert).
	\end{split}
	\end{equation*}
	As a result of the two above inequalities, as $v\to\infty$
	$$
	\tilde{\eta}_2(t)=O\left(\left\vert
	A\left(
	\frac{1}{1-F(x^*-1/y)}
	\right)\right\vert
	\right),
	$$
	Therefore, by regular variation of $1/(1-F(x^*-1/y))$ with index $-1/\gamma$, $\tilde{\eta}_2(y)$ is eventually dominated by a regularly varing function of index $-\rho/\gamma$. 
	The final result now follows. 
\end{proof}
}\fi
\begin{lem}\label{lem:densexpand2}
	If $\gamma<0$ and $\rho<0$, there exist $\tilde{\delta}>0$ such that, as $y\to \infty$,
	$$
	\frac{f(x^*-1/y)}{y^2}=(1-\gamma y)^{1/\gamma-1}\left[1+O(\{1-H_{-\gamma}(y)\}^{\tilde{\delta}})\right]
	$$
\end{lem}
\if1\blind
{
\begin{proof}
	The function $\tilde{f}(y):=f(x^*-1/y)y^{-2}$ is the density of the distribution function $\tilde{F}(y):=F(x^*-1/y)$, which is in the domain of attraction of $G_{\tilde{\gamma}}$, with $\tilde{\gamma}=-\gamma$. Moreover, 
	$$
	\tilde{\eta}(y)=\frac{(1+\tilde{\gamma} y)\tilde{f}(y)}{1-\tilde{F}(y)} -1.
	$$
	By Lemma \ref{lem:bigO2} and regular variation of $1-H_{\tilde{\gamma}}$ with index $-1/\tilde{\gamma}$, we have
	$$
	\tilde{\eta}(y)=O(\{1-H_{\tilde{\gamma}}(y)\}^{\tilde{\delta}})
	$$
	for any $\tilde{\delta}>0$ such that $-\tilde{\delta}/\tilde{\gamma}>\tilde{\varrho}$. Therefore, by  Proposition 2.1.5 in \cite{falk2010}, as $y\to\infty$ it holds that
	$$
	\tilde{f}(y)=h_{\tilde{\gamma}}(y)[1+O(\{1-H_{\tilde{\gamma}}(y)\}^{\tilde{\delta}})],
	$$
	which is the result.
\end{proof}
}\fi

{Finally, in} order to exploit Theorem \ref{theo:hellrate} to give bounds on Kullback-Leibler and higher order divergences, {we 
introduce by the next lemma
a uniform bound on density ratios. } 
%
	\begin{lem}\label{lem:boundens}
		Under the assumptions of Theorem \ref{theo:hellrate}, if $\rho <0$ and $\gamma \neq 0${,}
		then there exist a $t_1 <x^*$ and a constant $M \in (0,\infty)$ such that
		$$
		\sup_{t \geq t_1}\sup_{ 0<x<{\tilde{x}_t^*}}\frac{\tilde{l}_t(x)}{h_\gamma(x)}<M.
		$$
	\end{lem}
\if1\blind
{
\begin{proof}
We analyse the cases where $\gamma>0$ and $\gamma<0$ separately.

\textit{Case 1: $\gamma>0$}. In this case, 
{$\tilde{l}_t=l_t$.} By Lemma 4.4, there are positive constants $\kappa$, $\delta$ and $\epsilon$ such that, for all large $t$ and all $x>0$
\begin{equation*}
	\begin{split}
		\frac{l_t(x)}{h_\gamma(x)}&\leq \frac{h_{\gamma}(s(t)x+t)}{h_\gamma(x)}\frac{s(t)}{1-F(t)}
		\left[
		1+\kappa\left\lbrace
		1-H_\gamma(s(t)x+t)
		\right\rbrace^\delta
		\right]\\
		&\leq \left[
		\frac{1+\gamma x}{(1+\gamma t)/s(t) +\gamma x}
		\right]^{1+1/\gamma}\frac{1+\epsilon}{(s(t))^{1/\gamma}(1-F(t))}.
	\end{split}
\end{equation*}
Moreover, by Lemma 4.3 it holds that as $t\to\infty$
$$
\frac{1+\gamma t}{s(t)}=1+\eta(t)=1+o(1)
$$
and, in turn, $(s(t))^{1/\gamma}\sim (1+\gamma t)^{1/\gamma}$. These two facts, combined with the tail equivalence relation in formula \eqref{eq:tailequiv}, imply that for all sufficiently large $t$ and all $x>0$
\begin{equation*}
	\begin{split}
		\frac{l_t(x)}{h_\gamma(x)}&\leq\left[
		\frac{1+\gamma x}{ 1-\epsilon+\gamma x}
		\right]^{1+1/\gamma}\frac{1+\epsilon}{(1-\epsilon)\tau}\\
		&\leq \left[
		\frac{1}{ 1-\epsilon}
		\right]^{1+1/\gamma}\frac{1+\epsilon}{(1-\epsilon)\tau}.
	\end{split}
\end{equation*}
The result now follows.

\textit{Case 2: $\gamma<0$}. In this case, for any $x\in(0, 
{-1/\gamma}
)$
\begin{equation*}
	\begin{split}
		\tilde{l}_t(x)=f\left(x^*-\frac{1}{y}\right)\frac{1}{y^2}\frac{y^2 \tilde{s}(t)}{1-F(t)}
	\end{split}
\end{equation*}
where
$$
y\equiv y(x,t):={\frac{1}{{s}(t)}\left[-
	\frac{1}{\gamma}-x
	\right]^{-1}}
$$
Note that $y$ is bounded from below by 
{$-\gamma /s(t)$,}
which converges to $\infty$ as $t \to x^*$. Thus, by Lemma 4.6 there are positive constants $\tilde{\delta}$, $\epsilon$ and $\tilde{\kappa}$ such that
\begin{equation*}
	\begin{split}
		\tilde{l}_t(x)&\leq(1-\gamma y)^{1/\gamma-1}[1+\tilde{\kappa}\{1-H_{-\gamma}(y)\}^{\tilde{\delta}}]\frac{y^2 \tilde{s}(t)}{1-F(t)}\\
		&\leq { h_\gamma(x)}
		%
		\left[
		{s(t)\left(
			-\frac{1}{\gamma} -x
			\right)}
		-\gamma
		\right]^{\frac{1}{\gamma} -1}
		\frac{(1+\epsilon){(-\gamma^{-1} s(t))^{-1/\gamma}}}{1-F(t)}
		.
	\end{split}
\end{equation*}
By hypothesis, it holds that {$x<-1/\gamma$,}
%
thus
\begin{equation*}
	\begin{split}
		\left[
		{s(t)\left(
			-\frac{1}{\gamma} -x
			\right)}
		-\gamma
		\right]^{\frac{1}{\gamma} -1}
		\leq (-\gamma)^{\frac{1}{\gamma} -1}.
	\end{split}
\end{equation*}
Finally, for all large $t$,
$$
\frac{ -\gamma^{-1}s(t)}{x^*-t}\leq 
{(1+\epsilon)}
$$
Combining all the above inequalities we can now conclude that, for all large $t$ and for any $x\in(0, (x^*-t)/\tilde{s}(t))$,
\begin{equation*}
	\begin{split}
		\frac{\tilde{l}_t(x)}{h_\gamma(x)}&\leq (1+\epsilon)^{1-1/\gamma}(-\gamma)^{\frac{1}{\gamma}}\frac{(x^*-t)^{-\frac{1}{\gamma}}}{1-F(t)}.
	\end{split}
\end{equation*}
Now, setting $t=U(v)$, we have that $v\to\infty$  if and only if $t\to x^*$ and, by Theorem 2.3.6 in \cite{dehaan+f06}, there is a constant $\varpi>0$ such that for all large $t$
\begin{equation*}
	\begin{split}
		\frac{(x^*-t)^{-\frac{1}{\gamma}}}{1-F(t)}\leq v [(1+\epsilon)\varpi v^\gamma]^{-\frac{1}{\gamma}}=[(1+\epsilon)\varpi ]^{-\frac{1}{\gamma}}
	\end{split}
\end{equation*}
The result now follows.

\end{proof}
}\fi

	%
%
\if1\blind
{	
\begin{proof}

Note that for any $\gamma'>-1/2$ and $\sigma>0$
\begin{equation*}
\begin{split}
\hell(h_\gamma; h_{\gamma'}(\sigma \,\cdot \, ) \sigma)&\leq 
\sqrt{\int_{\mathbb{R}}\left[
	\nu_x(\gamma)-\nu_x(\gamma')
	\right]^2\diff x
}
+
\sqrt{\int_{\mathbb{R}}\left[
	\psi_{\gamma,x}(\sigma)-\psi_{\gamma,x}(1)
	\right]^2\diff x
}
\end{split}
\end{equation*}
In what follows, we bound the two terms on the right-hand side for $\gamma'\in (\gamma \pm \epsilon)$ and $\sigma\in (1\pm \epsilon)$, for a suitably small $\epsilon>0$. We study the the cases where $\gamma>0$, $\gamma<0$ and $\gamma=0$ separately.

\textit{Case 1: $\gamma>0$.} An application of the mean-value theorem and Lemma \ref{lem:auxlip1}\ref{res:pos} yields that, for a function $\xi(x)\in (\gamma\wedge\gamma',\gamma\vee\gamma')$,
\begin{equation*}
\begin{split}
\sqrt{\int_{\mathbb{R}}\left[
	\nu_x(\gamma)-\nu_x(\gamma')
	\right]^2\diff x}&=\vert \gamma - \gamma' \vert\sqrt{\int_0^\infty 
	\left[
	\nu_x'(\xi(x))
	\right]^2 \diff x}\\
&\leq 	\frac{\vert \gamma - \gamma'\vert}{2}
\sqrt{
	\int_0^\infty(1+x\xi(x))^{-3-\frac{1}{\xi(x)}}\left[
	\frac{x\ln(1+x\xi(x))}{\xi(x)}\right]^2\diff x
}\\
& 
\quad+
\frac{\vert \gamma - \gamma'\vert}{2}
\sqrt{
	\int_0^\infty(1+x\xi(x))^{-3-\frac{1}{\xi(x)}}x^2\diff x
}. 
\end{split}
\end{equation*}
On one hand, it holds that
\begin{equation*}
\begin{split}
&\int_0^\infty(1+x\xi(x))^{-3-\frac{1}{\xi(x)}}\left[
\frac{x\ln(1+x\xi(x))}{\xi(x)}\right]^2\diff x\\
&\leq 4\int_0^\infty(1+x(\gamma-\epsilon))^{-1-\frac{1}{\gamma+\epsilon}}\left[
\frac{\ln(1+x(\gamma-\epsilon))}{(\gamma-\epsilon)^2}\right]^2\diff x\\
&\leq 8\frac{(\gamma+\epsilon)^3}{(\gamma-\epsilon)^5}. 
\end{split}
\end{equation*}
On the other hand, it holds that
\begin{equation*}
\begin{split}
&\int_0^\infty(1+\xi(x))^{-3-\frac{1}{\xi(x)}}x^2\diff x\\
&\leq \int_0^\infty(1+x(\gamma-\epsilon))^{-1-\frac{1}{\gamma+\epsilon}}\frac{1}{(\gamma-\epsilon)^2}\diff x\\
&\leq \frac{(\gamma+\epsilon)}{(\gamma-\epsilon)^3}. 
\end{split}
\end{equation*}
While, an application of the mean-value theorem and Lemma \ref{lem:auxlipbis}\ref{lem:auxlip2} yields that, for a function $\varsigma(x)\in(1\wedge \sigma,1\vee \sigma)$,
\begin{equation*}
\begin{split}
\int_{\mathbb{R}}\left[
\psi_{\gamma',x}(\sigma)-\psi_{\gamma,x}(1)
\right]^2\diff x&=	
\int_0^\infty \left[
\psi_{\gamma,x}'(\varsigma(x))
\right]^2\diff x\\
&	\leq \left(
\frac{1}{\gamma^2}+\sqrt{\frac{1}{2\gamma +1}}
\right)^2\left(\frac{1+\epsilon}{1-\epsilon}\right)^{5}.
\end{split}
\end{equation*}
The result now follows.

\textit{Case 2: $\gamma<0$.} Assume that $\gamma<\gamma'$, then an application of the mean-value theorem and Lemma \ref{lem:auxlip1}\ref{res:neg} yields that, for a function $\xi(x)\in (\gamma,\gamma')$,
\begin{equation*}
\begin{split}
\int_{\mathbb{R}}\left[
\nu_x(\gamma)-\nu_x(\gamma')
\right]^2\diff x&=\vert \gamma - \gamma' \vert^2\int_0^{-1/\gamma} 
\left[
\nu_x'(\xi(x))
\right]^2 \diff x+1-H_{\gamma'}(-1/\gamma)\\
&\leq 	\frac{\vert \gamma - \gamma'\vert^2}{4}\left[
\sqrt{
	\int_0^{-1/\gamma}(1+x\xi(x))^{-3-\frac{1}{\xi(x)}}x^4\diff x
}\right.\\
& 
\quad+\left.
\sqrt{
	\int_0^{-1/\gamma}(1+x\xi(x))^{-3-\frac{1}{\xi(x)}}x^2\diff x
}\right]^2 + 1-H_{\gamma'}(-1/\gamma).
\end{split}
\end{equation*}
First, for a constant $\beta$ satisfying $0< \beta < 1/(\epsilon-\gamma)-2$, we have that
\begin{equation*}
\begin{split}
\int_0^{-1/\gamma}(1+x\xi(x))^{-3-\frac{1}{\xi(x)}}x^4\diff x&\leq\frac{1}{\gamma^4}\int_0^{-1/\gamma}(1+\gamma x)^{-1+\beta}\diff x\\
&\leq \frac{1}{(-\gamma)^5}\frac{1}{\beta}.
\end{split}
\end{equation*}
Similarly, 
\begin{equation*}
\begin{split}
\int_0^{-1/\gamma}(1+x\xi(x))^{-3-\frac{1}{\xi(x)}}x^2\diff x&\leq\frac{1}{\gamma^2}\int_0^{-1/\gamma}(1+\gamma x)^{-1+\beta}\diff x\\
&\leq \frac{1}{(-\gamma)^3}\frac{1}{\beta}.
\end{split}
\end{equation*}
Finally, if $\epsilon$ is small enough, $ 1-H_{\gamma'}(-1/\gamma)\leq (1-\gamma'/\gamma)^2$. Thus, we can conclude that
\begin{equation*}
\begin{split}
\int_{\mathbb{R}}\left[
\nu_x(\gamma)-\nu_x(\gamma')
\right]^2\diff x&\leq \vert \gamma - \gamma'\vert^2 \frac{1+1/2\beta}{(-\gamma)^5}.
\end{split}
\end{equation*}
A similar reasoning when $\gamma >\gamma'$ yields that
\begin{equation*}
\begin{split}
\int_{\mathbb{R}}\left[
\nu_x(\gamma)-\nu_x(\gamma')
\right]^2\diff x&\leq \vert \gamma - \gamma'\vert^2 \frac{1+2/\beta}{(-\gamma')^5}\\
&\leq 
\vert \gamma - \gamma'\vert^2 \frac{1+2/\beta}{(-\gamma-\epsilon)^5}
.
\end{split}
\end{equation*}

Next, assuming that $\sigma <1$,  an application of the mean-value theorem and the first half of Lemma \ref{lem:auxlipbis}\ref{lem:auxlip3} yields that for a function $\varsigma(x  )\in (1-\epsilon,1)$ and a constant $\zeta>0$
\begin{equation*}
\begin{split}
\int_{\mathbb{R}}\left[
\psi_{\gamma,x}(\sigma)-\psi_{\gamma,x}(1)
\right]^2\diff x&=
(1-\sigma)^2\int_0^{-\sigma/\gamma} \left[
\psi_{\gamma,x}'(\varsigma(x))
\right]^2\diff x+(1-\sigma)^{-1/\gamma}\\
&\leq (1-\sigma)^2\left[\frac{1}{\zeta (1-\epsilon)^2(-\gamma)}\left(
\frac{1}{\gamma^2}+1
\right)^2
+1
\right].
\end{split}
\end{equation*}
While, if $\sigma>1$, for a function $\varsigma(x  )\in (1,1+\epsilon)$
\begin{equation*}
\begin{split}
\int_{\mathbb{R}}\left[
\psi_{\gamma,x}(\sigma)-\psi_{\gamma,x}(1)
\right]^2\diff x&=
(1-\sigma)^2\int_0^{-1/\gamma} \left[
\psi_{\gamma,x}'(\varsigma(x))
\right]^2\diff x+(1-1/\sigma)^{-1/\gamma}\\
&\leq (1-\sigma)^2\left[
\frac{1}{\zeta (-\gamma)}\left(
\frac{1}{\gamma^2}+1
\right)^2
+1
\right].
\end{split}
\end{equation*}
The result now follows.

\textit{Case 3: $\gamma=0$.}  Assume that $\gamma'>0$, then an application of the mean-value theorem and Lemma \ref{lem:auxlip1}\ref{res:pos} yields that, for a function $\xi(x)\in (0,\gamma')$,
\begin{equation*}
\begin{split}
\int_{\mathbb{R}}\left[
\nu_x(\gamma)-\nu_x(\gamma')
\right]^2\diff x&=\vert \gamma - \gamma' \vert^2\int_0^{\infty} 
\left[
\nu_x'(\xi(x))
\right]^2 \diff x\\
&\leq 	\frac{\vert \gamma - \gamma'\vert^2}{4}\left[
\sqrt{
	\int_0^{\infty}
	(1+x\xi(x))^{-3-\frac{1}{\xi(x)}}
	\left[
	\frac{x\ln(1+x\xi(x))}{\xi(x)}\right]^2\diff x
}\right.\\
& 
\quad+\left.
\sqrt{
	\int_0^{-\infty}(1+x\xi(x))^{-3-\frac{1}{\xi(x)}}x^2\diff x
}\right]^2.
\end{split}
\end{equation*}
On one hand, we have
\begin{equation*}
\begin{split}
\int_0^{\infty}
(1+x\xi(x))^{-3-\frac{1}{\xi(x)}}
\left[
\frac{x\ln(1+x\xi(x))}{\xi(x)}\right]^2\diff x &\leq \int_0^{\infty}
(1+x\xi(x))^{-3-\frac{1}{\xi(x)}}
x^4\diff x\\
&\leq  \int_0^{\infty}
(1+x\gamma')^{-3-\frac{1}{\gamma'}}
x^4\diff x+\int_0^\infty e^{-x}x^4\diff x\\
&\leq 36 + \Gamma(5).
\end{split}
\end{equation*}
On the other hand, we have
\begin{equation*}
\begin{split}
\int_0^{\infty}
(1+x\xi(x))^{-3-\frac{1}{\xi(x)}}
x^2\diff x
&\leq  \int_0^{\infty}
(1+x\gamma')^{-3-\frac{1}{\gamma'}}
x^2\diff x+\int_0^\infty e^{-x}x^2\diff x\\
&\leq 3 + \Gamma(3).
\end{split}
\end{equation*}

Assume next that $\gamma'<0$, then an application of the mean-value theorem and Lemma \ref{lem:auxlip1}\ref{res:neg} yields that, for a function $\xi(x)\in (-\epsilon,0)$,
\begin{equation*}
\begin{split}
\int_{\mathbb{R}}\left[
\nu_x(\gamma)-\nu_x(\gamma')
\right]^2\diff x&=\vert \gamma - \gamma' \vert^2\int_0^{-1/\gamma'} 
\left[
\nu_x'(\xi(x))
\right]^2 \diff x+e^{1/\gamma'}\\
&\leq 	\frac{\vert \gamma - \gamma'\vert^2}{4}\left[
\sqrt{
	\int_0^{-1/\gamma'}(1+x\xi(x))^{-3-\frac{1}{\xi(x)}}x^4\diff x
}\right.\\
& 
\quad+\left.
\sqrt{
	\int_0^{-1/\gamma'}(1+x\xi(x))^{-3-\frac{1}{\xi(x)}}x^2\diff x
}\right]^2 + e^{1/\gamma'}.
\end{split}
\end{equation*}
On one hand, for $\epsilon$ sufficiently small we have
\begin{equation*}
\begin{split}
\int_0^{-1/\gamma'}(1+x\xi(x))^{-3-\frac{1}{\xi(x)}}x^4\diff x &\leq 
\int_0^{-1/\gamma'}
(1+x\gamma')^{-3-\frac{1}{\gamma'}}
x^4\diff x+\int_0^{-1/\gamma'} e^{-x}x^4\diff x\\
&\leq \frac{13}{2}\Gamma(5)  
\end{split}
\end{equation*}
and
\begin{equation*}
\begin{split}
\int_0^{-1/\gamma'}(1+x\xi(x))^{-3-\frac{1}{\xi(x)}}x^2\diff x &\leq 
\int_0^{-1/\gamma'}
(1+x\gamma')^{-3-\frac{1}{\gamma'}}
x^4\diff x+\int_0^{-1/\gamma'} e^{-x}x^2\diff x\\
&\leq \frac{3}{2}\Gamma(3).
\end{split}
\end{equation*}
On the other hand, for $\epsilon$ sufficiently small we have $e^{1/\gamma'}\leq \vert\gamma'-\gamma\vert^2$.

Finally, some algebraic manipulations yield
\begin{equation*}
\begin{split}
\int_{\mathbb{R}}\left[
\psi_{\gamma,x}(\sigma)-\psi_{\gamma,x}(1)
\right]^2\diff x&=\int_0^\infty\left[
\sqrt{e^{-x/\sigma}\frac{1}{\sigma}}
-\sqrt{e^{-x}}
\right]^2\diff x
\\
&\leq
\frac{(1-\sigma)^2}{(1-\epsilon)^2}
\left[1+\frac{1}{2}\left(\frac{1+\epsilon}{1-\epsilon}\right)^{3/2}
\right]^2.
\end{split}
\end{equation*}
The proof is now complete.

\end{proof}
}\fi

\subsection{Proof of Theorem \ref{theo:hellrate}}
For every $x_t>0$, it holds that
\begin{equation*}
\begin{split}
\hell^2(l_t; h_\gamma)&=\int_0^{x_t}+\int_{x_t}^\infty
\left[
\sqrt{f_t}(x)-\sqrt{h_\gamma({(x-c(t))}/s(t))
	/s(t)
}
\right]^2 \diff x\\
&\leq \int_{0}^{\phi_t(x_t)}e^{-y}
\left[
1-\sqrt{
	e^{{-}q_t(y)}(1+q_t'(y))
}
\right]^2\diff y\\
&\quad + 
\left[
\sqrt{1-F_t(x_t)}
+	\sqrt{1-{H_\gamma\left(
	\frac{x_t-c(t)}{s(t)}
	\right)}
}
\right]^2\\
&=:\mathcal{I}_1(t)+\mathcal{I}_2(t){.}
\end{split}
\end{equation*}
Let $x_t$ be such that the following equality holds 
$$
\phi_t(x_t)=-\alpha \ln \vert
A(e^{T^{-1}(t)})
\vert,
$$
for a positive constant $\alpha$ to be specified later. 
Then, by Lemmas \ref{lem:aux1}-\ref{lem:aux2}, for a suitably small $\varepsilon>0$ {there exist $\kappa_3>0$}
such that for all sufficiently large $t$
\begin{equation*}
\begin{split}
\mathcal{I}_1(t)&\leq 
\int_0^{-\alpha \ln \vert
	A(e^{T^{-1}(t)})\vert}\kappa_3
\vert
A(e^{T^{-1}(t)})
\vert^2 e^{(4\varepsilon-1) y}\diff y
\\
&\leq 
\kappa_3
\vert
A(e^{T^{-1}(t)})
\vert^2 \left[
1-\vert
A(e^{T^{-1}(t)})
\vert^{\alpha_1}
\right]
\end{split}
\end{equation*}
where $\alpha_1:=\alpha(1-4\varepsilon)$ 
{is} positive. Moreover, on one hand we have the identity 
$$
1-F_t(x_t)=\vert
A(e^{T^{-1}(t)})
\vert^\alpha.
$$
On the other hand, for some constan{t $\kappa_5>0$} we have the inequality
\begin{equation*}
\begin{split}
1-
{H_\gamma\left(
	\frac{x_t-c(t)}{s(t)}
	\right)}&=\vert
A(e^{T^{-1}(t)})
\vert^\alpha\exp\left\lbrace
-q_t\left(-\alpha \ln \vert
A(e^{T^{-1}(t)})\vert \right)
\right\rbrace\\
&\leq
\vert
A(e^{T^{-1}(t)})
\vert^\alpha \exp
\left\lbrace\kappa_5 	\vert
A(e^{T^{-1}(t)})
\vert^{1-2\varepsilon\alpha} \right\rbrace{.}
\end{split}
\end{equation*}
Consequently,
$$
\mathcal{I}_2(t)\leq
\vert A(e^{T^{-1}(t)})
\vert^\alpha \left[
1+ \exp
\left\lbrace\frac{\kappa_5}{2} 	\vert
A(e^{T^{-1}(t)})
\vert^{1-2\varepsilon\alpha} \right\rbrace
\right]
.
$$
Now,
 we can choose  $\alpha >2$ and $\varepsilon$ small enough, so that 
$$	\vert
A(e^{T^{-1}(t)})
\vert^\alpha < 	\vert
A(e^{T^{-1}(t)})
\vert^2 
$$
and  $\alpha_2:=1-2\varepsilon \alpha>0$. 
The conclusion then follows noting that
$
T^{-1}(t)=-\ln(1-F(t)) 
$
and, in turn,
$$
\vert
A(e^{T^{-1}(t)})
\vert=\vert A(v)\vert.
$$

\subsection{Proof of Proposition \ref{prop:boundGPD}}
Assume first that for all large $t$
\begin{equation*}\label{eq:mu}
	\mu_t=\frac{x^*-t}{s(t)}+\frac{1}{\gamma}
\end{equation*}
is positive. In this case, we have the following identities
\begin{equation*}
	\begin{split}
		\hell^2\left(
		h_\gamma, h_\gamma\left(
		\, \cdot \, - \mu_t
		\right)
		\right) &= \int_0^{-1/\gamma}
		\left[
		\sqrt{h_\gamma(x)}-\sqrt{h_\gamma(x-\mu_t)}
		\right]^2\diff x+1-H_\gamma \left( 
		-\frac{1}{\gamma}-\mu_t
		\right)\\
		&= \int_0^\infty e^{-s} \left[
		\sqrt{\frac{e^{-\frac{1}{\gamma} \ln
		\left(
		1- \gamma e^{-\gamma s} \mu_t
		\right)	
		}}{1-\gamma e^{-\gamma s} \mu_t}
		}
		-1
		\right]^2 \diff s+ (-\gamma \mu_t)^{-1/\gamma}.
	\end{split}
\end{equation*}
Concerning the first term on the right-hand side, for all $s>0$ as $t\to x^*$ we have
$$
0 \leq \ln 
\left(
1- \gamma e^{-\gamma s} \mu_t
\right)	 \leq \ln
\left(
1- \gamma \mu_t
\right)	= O(\mu_t)
$$
and
$$
1 \geq \frac{1}{1- \gamma e^{-\gamma s} \mu_t} \geq  \frac{1}{1- \gamma  \mu_t}=1+O(\mu_t),
$$
where, for a positive constant $\tau$, 
$
{\mu_t}
$
satisfies
$$
\frac{\mu_t}{|A(v)|} =(1+o(1))\tau
$$
 \cite[e.g.,][Lemma 4.5.4]{dehaan+f06}. Then, as $t \to x^*$
\begin{equation*}
	\begin{split}
		\int_0^\infty e^{-s} \left[
		\sqrt{\frac{e^{-\frac{1}{\gamma} \ln
					\left(
					1- \gamma e^{-\gamma s} \mu_t
					\right)	
			}}{1-\gamma e^{-\gamma s} \mu_t}
		}
		-1
		\right]^2 \diff s=O(\mu_t^2)=O(|A(v)|^2).
	\end{split}
\end{equation*}
Concerning the second term, as $t \to x^*$ we have
$$
(-\gamma \mu_t)^{-1\gamma}= (-\gamma)^{-1/\gamma} (1+o(t))\tau^{-1/\gamma} |A(v)|^{-1/\gamma}.
$$
The result now follows for the case where $\mu_t$ is ultimately positive. When it is ultimately negative, simply note that
\begin{equation}\label{eq:hell_id}
\hell^2\left(
h_\gamma, h_\gamma\left(
\, \cdot \, - \mu_t
\right)
\right) =
\hell^2\left(
h_\gamma\left(
\, \cdot \,  +\mu_t
\right),
h_\gamma
\right) 
=\hell^2\left(
h_\gamma\left(
\, \cdot \, -(- \mu_t)
\right),
h_\gamma
\right) 
\end{equation}
and proceed as above, but replacing $\mu_t$ with $-\mu_t$.

\subsection{Proof of Corollary \ref{cor:rate}}
Observe that, for $\gamma <0$,
$$
\hell(l_t, h_\gamma)=\hell(\tilde{l}_t, h_\gamma(\, \cdot\,  +\mu_t)).
$$
Moreover, note that by triangular inequality and the first identity in \eqref{eq:hell_id}
\begin{equation*}
	\begin{split}
		&\hell\left(
		h_\gamma, h_\gamma\left(
		\, \cdot \, - \mu_t
		\right)
		\right)- \hell(\tilde{l}_t, h_\gamma(\, \cdot\,  +\mu_t))\\
		&\quad \leq \hell(\tilde{l}_t, h_\gamma(\, \cdot\,  +\mu_t))\\
		& \quad \leq \hell(	h_\gamma, h_\gamma\left(
		\, \cdot \, - \mu_t
		\right))+ \hell(\tilde{l}_t, h_\gamma(\, \cdot\,  +\mu_t)).
	\end{split}
\end{equation*}
Whenever $\gamma \geq -1/2$, by Theorem \ref{theo:hellrate} and Proposition \ref{prop:boundGPD}, the term on the third line is of order $O(|A(v)|)$ as $t \to x^*$. Instead, if $\gamma <-1/2$, by Theorem \ref{theo:hellrate} and Proposition \ref{prop:boundGPD}, as $t\to x^*$ the term on the first line satisfies
\begin{equation*}
	\begin{split}
		\hell\left(
		h_\gamma, h_\gamma\left(
		\, \cdot \, - \mu_t
		\right)
		\right)- \hell(\tilde{l}_t, h_\gamma(\, \cdot\,  +\mu_t)) &\geq c_3 |A(v)|^{-1/2\gamma} - O(|A(v)|)\\
		&=c_3 |A(v)|^{-1/2\gamma}(1+o(1)).
	\end{split}
\end{equation*}
The result now follows.

\subsection{Proof of Corollary \ref{cor:kulbrate}}
By Lemma 8.2 in \cite{10.1214/aos/1016218228} 
$$
\kulb(\tilde{l}_t; h_\gamma)\leq 2 \left[
\sup_{ 0<x<
{\tilde{x}_t^*}}\frac{\tilde{l}_t(x)}{h_\gamma(x)}
\right]\hell^2(\tilde{l}_t; h_\gamma).
$$
Moreover, by Lemma B.3 in \cite{ghosal2017}, for $p\geq 2$
$$
\pdiv_p(\tilde{l}_t; h_\gamma)\leq 2 p! \left[ \sup_{ 0<x<{\tilde{x}_t^*}}\frac{\tilde{l}_t(x)}{h_\gamma(x)}
\right]\hell^2(\tilde{l}_t; h_\gamma).
$$
The conclusion now follows by combining the above inequalities and applying Theorem \ref{theo:hellrate} and Lemma \ref{lem:boundens}.
\section*{Acknowledgements}
Simone Padoan is supported by the Bocconi Institute for Data Science and Analytics (BIDSA), Italy. 
 
\bibliographystyle{chicago} 
\bibliography{bibliopm_final3}

\end{document}